\documentclass[journal,transmag]{IEEEtran}

\usepackage{graphicx}
\usepackage{amsmath}
\usepackage{tuda-pgfplots}
\usepackage{tablefootnote}
\usepackage{siunitx}
\usepackage[font=footnotesize]{caption}
\usepackage{subcaption}
\usepackage{amsfonts}
\usepackage{tikz}
\usepackage{pgfplots}
\usepackage{cite}
\usepackage{xcolor}
\usepackage{comment}
\usepackage{mathtools}
\usepackage{hyperref}
\usepackage{eso-pic}
\hypersetup{hidelinks}
\pgfplotsset{compat=1.18}

\graphicspath{ {./images/} }
\DeclareGraphicsExtensions{.pdf,.jpeg,.png,.tikz}

\newcommand{\CURL}{\ensuremath{\operatorname{\iflanguage{ngerman}{rot}{curl}}}}
\newcommand{\DIV}{\ensuremath{\operatorname{div}}}
\newcommand{\Hilbert}{\ensuremath{\mathrm{H}}}
\newcommand{\HSpace}[3]{\ensuremath{\Hilbert^{#1}_{#2}\!\left(#3\right)}}
\newcommand{\HSPACE}[4]{\HSpace{#1}{#2}{#3,#4}}
\newcommand{\rv}      {\ensuremath{\vec{r}}} \newcommand{\Afield}  {\ensuremath{\vec{A}}}   \newcommand{\Jfield}  {\ensuremath{\vec{J}}}   \newcommand{\Vo}{\ensuremath{V_{\mathrm{ext}}}} \newcommand{\Vi}{\ensuremath{V_{\mathrm{int}}}} \newcommand{\Gio}{\ensuremath{\Gamma}}\newcommand{\As} {\ensuremath{\Afield_\mathrm{s}}} \newcommand{\Am}{\ensuremath{\Afield_\mathrm{m}}} \newcommand{\Ag}{\ensuremath{\Afield_\mathrm{g}}} \newcommand{\Agh}{\ensuremath{\Afield_{\mathrm{g},h}}} \newcommand{\Js} {\ensuremath{\Jfield_\mathrm{s}}} \newcommand{\Kg}{\ensuremath{\vec{K}_\mathrm{g}}} \newcommand{\Kgh}{\ensuremath{\vec{K}_{\mathrm{g},h}}} 

\makeatletter
\patchcmd{\@makecaption}
{\scshape}
{}
{}
{}
\makeatother

\AddToShipoutPicture*{
	\footnotesize\sffamily\raisebox{0.8cm}{\hspace{1.4cm}\fbox{
			\parbox{\textwidth}{
				© 2026 IEEE. Personal use of this material is permitted. Permission from IEEE must be obtained for all other uses, including reprinting/republishing this material for advertising or promotional purposes, collecting new collected works for resale or redistribution to servers or lists, or reuse of any copyrighted component of this work in other works.
			}
	}}
}
\begin{document}
\title{A Reduced Magnetic Vector Potential Approach with Higher-Order
Splines}

\author{\IEEEauthorblockN{Merle Backmeyer\IEEEauthorrefmark{1,2},
Laura A. M. D'Angelo\IEEEauthorrefmark{1}, Brahim Ramdane\IEEEauthorrefmark{2} and
Sebastian Schöps\IEEEauthorrefmark{1}}\IEEEauthorblockA{\IEEEauthorrefmark{1}Computational Electromagnetics Group, Technische Universität Darmstadt, 64289 Darmstadt, Germany \\
\IEEEauthorrefmark{2} Univ. Grenoble Alpes, CNRS, Grenoble INP\footnote{Institute of Engineering Univ. Grenoble Alpes}, G2Elab, 38000 Grenoble, France}\thanks{Received 25 February 2026; revised 11 June 2026 ; accepted 13 July 2026 . Date of publication 22 July 2026; date of current version July 2026. Corresponding author: M. Backmeyer (email: merle.backmeyer@univ-grenoble-alpes.fr).}
\thanks{Color versions of one or more figures in this article are available at
https://doi.org/10.1109/TMAG.2026.3716008.}
\thanks{Digital Object Identifier 10.1109/TMAG.2026.3716008}
}

\markboth{A Reduced Magnetic Vector Potential Approach with Higher-Order Splines}{Backmeyer \MakeLowercase{\textit{et al.}}: A Reduced Magnetic Vector Potential Approach with Higher-Order Splines}

\IEEEtitleabstractindextext{\begin{abstract}
This work presents a high‑order isogeometric formulation for magnetoquasistatic eddy‑current problems based on a decomposition into Biot–Savart–driven source fields and finite‑element reaction fields. Building upon a recently proposed surface‑only Biot–Savart evaluation, we generalize the reduced magnetic vector potential framework to the quasistatic regime and introduce a consistent high‑order spline discretization. The resulting method avoids coil meshing, supports arbitrary winding paths, and enables high‑order field approximation within a reduced computational domain. Beyond establishing optimal convergence rates, the numerical investigation identifies the requirements necessary to recover high‑order accuracy in practice, including geometric regularity of the enclosing interface, accurate kernel quadrature, and compatible trace spaces for the source–reaction coupling.
\end{abstract}

}

\maketitle
\IEEEdisplaynontitleabstractindextext
\IEEEpeerreviewmaketitle

\section{Introduction}
\label{sec:introduction}
\IEEEPARstart{W}{hen} solving eddy current problems, resolving filamentary coils in the computational mesh can be computationally inconvenient due to their high-aspect ratio and therefore the demand for very fine mesh resolution in their proximity. It is particularly challenging in isogeometric analysis (IGA), where conforming multi-patch discretizations inherit the tensor-product structure of NURBS and make thin features in the mesh particularly challenging to represent and refine efficiently~\cite{Hughes_2005aa}.
To avoid modeling individual turns, homogenization and coil surrogate models (e.g., solid or stranded conductor models) have long been employed and continue to be refined \cite{Dyck_2004aa, Sabariego_2017aa, Paakkunainen_2024ad}. These approaches are especially attractive when the coil itself is the region of interest and one needs to compute current distributions within the winding pack. When the coil merely acts as a source through its far field, a different strategy is to exclude the coil from the finite-element domain. Instead, the magnetic field is decomposed into a prescribed linear source field and a potentially non-linear reaction field induced in the conducting region. The source field can be computed via a fundamental solution, the Biot-Savart law, while the reaction field is computed on a mesh that resolves only the conductors and air region, not the coil geometry itself. This “non-resolved coil” strategy is known as the reduced magnetic vector potential formulation (RMVP) was originally proposed in \cite{Biro_1999aa}. It is worth noting that this field‑decomposition viewpoint, in which the source‑generated field in a linear subdomain is represented using the fundamental solution, is not unique to RMVP. Closely related ideas appear in volume‑integral equations, boundary‑integral formulations, and in FEM–BEM coupling, where the exterior (or source) field is represented via the fundamental solution while only the reaction field is discretized numerically~\cite{Kress_1989aa, Albanese_1988aa, Hiptmair_2003aa}. These connections place RMVP within a broader class of classical electromagnetic formulations that exploit Green’s functions to reduce the computational domain and improve efficiency. However, this comes with its own challenge: the expensive Biot–Savart integrals have to be evaluated for each quadrature point of each element of the mesh. This can be mitigated by compression, e.g. fast multipole methods (FMM)~\cite{Greengard_1987aa}, or reducing the domain over which the fundamental solution is evaluated. A recent idea, proposed in \cite{DAngelo_2024ab} is to restrict the evaluation of the fundamental solution to a closed surface that separates the linear air domain from the magnetic region which typically features a nonlinear magnetic material behavior. We will refer to it as the interface RMVP. This surface-based representation preserves the correct source field in the exterior domain while drastically reducing the number of kernel evaluations. While it has demonstrated efficiency for magnetostatic problems with low‑order discretizations, the reported convergence remained suboptimal, leaving open questions regarding its behavior for higher-order discretizations and more general settings. Building on this foundation, we address these limitations. In particular, we show that optimal convergence is in fact attainable and identify the conditions required to achieve it. In addition, we show that the formulation can be directly applied to the magnetoquasistatic (MQS) regime, supports arbitrary winding paths and naturally accomodates higher-order basis functions. Higher-order discretizations promise reduced error per degree of freedom and better accuracy for field derivatives — critical for loss estimation, local flux metrics, and design sensitivities~\cite{Evans_2009aa}. However, to actually achieve optimal convergence, the entire chain must be high-order consistent: geometry approximation of the enclosing surface, high-order accurate kernel quadrature and compatible trace spaces for coupling that preserve the discretization’s asymptotic rates. This paper makes these requirements explicit and provides a step-by-step pathway to restore optimal convergence in practice. The implementation is based on open-source software and the code is publicly available \cite{Backmeyer_2026aa}. The remainder of this paper is structured as follows: In \autoref{sec:methodology}, the interface RMVP is revisited, higher-order discretization is introduced and the expected optimal convergence rates are discussed. \autoref{sec:results}
presents numerical results and \autoref{sec:conclusion} concludes the work. 

\section{Methodology}
\label{sec:methodology}

In this chapter, we revisit the methodology originally developed in~\cite{DAngelo_2024ab} within the framework of linear finite elements. As an extension of previous work, we develop a higher‑order isogeometric discretization and provide convergence estimates. In practice, the approach requires (i) robust construction of the enclosing surface; (ii) accurate numerical quadrature for the kernel interactions; and (iii) a consistent transfer of the source data onto the reaction problem. These steps become increasingly delicate for higher-order discretizations, where inadequate quadrature or function spaces can compromise the expected benefits of high-order approximations. We will visit them carefully in the following.

\subsection{Interface Reduced Magnetic Vector Potential Approach}
\begin{figure*}[!t]
\centering
\subfloat[Cross-section in the $xy$-plane]{\begin{tikzpicture}[scale=0.9]
				\draw[draw=black, fill=black, fill opacity=0.1, line width=1pt] (0,0) circle (2.5cm);
                \draw[fill=white,dashed] (0,0) circle (1.35cm);
				\draw[fill=TUDa-2c, opacity=0.3] (0,0) circle (1.2cm);
				\node at (-0.5,1.0) {$\Gio$};
				\node at (0,0) {$\Vi$}; \node at (2.1,0) {$\Vo$}; \draw[draw=TUDa-9b, fill=none, line width = 1pt] (0,0) circle (1.7cm);
\end{tikzpicture} \label{fig:xy-cut}}
\hfil
\subfloat[Cross-section in the $xz$-plane]{\begin{tikzpicture}[scale=0.9]
				\draw[draw=black, fill =black, fill opacity=0.1, line width=1pt] (-2cm,-2.7cm) rectangle ++(4cm, 5.4cm);
                \draw[fill=white, dashed] (-1cm,-1.5cm) rectangle ++(2cm,3cm);
				\draw[fill=TUDa-2c, opacity=0.3] (-0.5cm,-1cm) rectangle ++(1cm,2cm);
				\node at (-1.2cm,1.2cm) {$\Gio$};
				\node at (-0.4cm,0.75cm) {$\Vi$}; \node at (1.6,0.75cm) {$\Vo$};

\tikzset{coilpt/.style={fill=red, draw=red}}
\def\pts{1pt}
            
\coordinate (LTop) at (-1.35,  0.25);
              \coordinate (LBot) at (-1.35,  -0.25);
            
\foreach \P in {(-1.35,  0.25), (-1.35,  0.15),
                              (-1.35,  0.05),
                              (-1.35,  -0.05),
                              (-1.35,  -0.15),
                              (-1.35,  -0.25)} {\fill[coilpt] \P circle (\pts);}
            
\foreach \P in {( 1.35,  0.2),
                              ( 1.35,  0.1),
                              ( 1.35,  0),
                              ( 1.35, -0.1),
                              ( 1.35, -0.2)}
                {\fill[coilpt] \P circle (\pts);}
            
\coordinate (LLeftTop) at (-1.45,  0.25);
              \coordinate (LLeftBot) at (-1.45, -0.25);
            
              \draw[red, thick, line cap=round]
                (LLeftTop) -- (LTop)      (LLeftTop) -- (LLeftBot)  (LLeftBot) -- (LBot);     \end{tikzpicture}
 \label{fig:xz-cut}}
\caption{Problem setup illustrated through cross‑sections in the $xy$‑ and $xz$‑planes.  The closed surface $\Gio$ partitions the computational domain $V$ into the outer region $\Vo$ and the inner region $\Vi$. The conductor (blue) is fully contained in $\Vi$, while the coil (red), shown here with five turns, carries the source current in $\Vo$. This setup is analyzed in \autoref{sec:results} with two different coil configurations.}
\label{fig:problem_sketch}
\end{figure*}
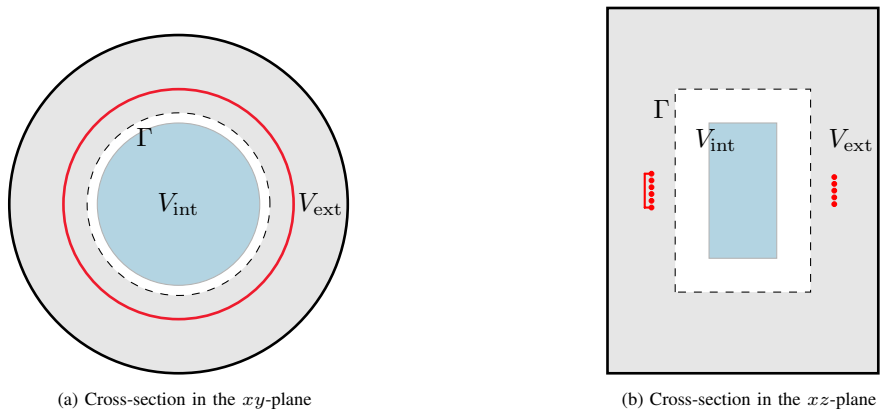

Let us assume that the computational domain $V = \Vo \cup \Vi \subset \mathbb{R}^3$ is composed of an outer region $\Vo$ that contains the coil embedded in a linear homogeneous material and an inner region $\Vi$ that contains the active regions (conductive with electrical conductivity $\sigma$ and/or magnetic with permeability $\mu$). The two domains are separated by the surface $\Gio$. This setup is depicted in~\autoref{fig:problem_sketch}.

We operate in the MQS regime, where displacement currents are disregarded. For clarity of presentation, we restrict our attention to the time‑harmonic regime, assuming linear material behavior also in $\Vi$ and sinusoidal source current density. This simplification is only introduced for conciseness: nonlinear and time‑dependent formulations can be handled using standard time‑stepping schemes, single or multistep. 

Before introducing the formal formulation, we briefly outline the underlying idea. Resolving the coil geometry inside the computational domain is avoided by representing the source term not in terms of the source current density $\Js$, but the source field $\As$ in $V$ that can be computed via the law of Biot-Savart in a preprocessing step. That is the idea of RMVP. This approach is valid, when the far field of the coil is of interest. In the interface RMVP approach, this source field is represented through an equivalent surface current density on an interface $\Gamma$ that separates the air region containing the source coils from the conducting region. To this end, the total field is decomposed such that the inner region is shielded from the direct influence of the source field. This is achieved by introducing an auxiliary image field $\Am$ in the outer region $\Vo$ such that the superposition $\As + \Am$ satisfies a homogeneous Dirichlet condition on the interface $\Gamma$. This procedure induces a virtual surface current density $\Kg$ on $\Gamma$. Although $\Kg$ does not correspond to a physical current on the interface, it contains all information required to reproduce the effect of the source field inside $\Vi$ based on the equivalence principle \cite[pp. 106-108]{Harrington_2001aa}. For more details, the reader is referred to \cite{DAngelo_2024ab}.

The magnetic vector potential $\Afield$ is decomposed into 
\begin{equation}
	\Afield = \begin{cases} \As + \Am + \Ag & \text{in } \Vo. \\
		\Ag & \text{in } \Vi, \end{cases}
        \label{eq:field_superpos}
\end{equation} where 
a source field $\As$ is induced by the coil currents, an image field $\Am$ that in superposition with the source field gives a fundamental solution in $\Vo$ that respects the homogeneous Dirichlet conditions at~$\Gio$, and a reaction field $\Ag$ that satisfies the eddy current equations in the computational domain with material constitutive laws (possibly nonlinear), driven by a surface current that encodes the source contribution. We compute the source field for any (quadrature) point on the surface $\Gio$ (and any point of interest in~$\Vo$) via the Biot–Savart kernel~\cite[Sec.~5.4]{Jackson_1998aa}, i.e.
\begin{equation}
		\As(\rv) = \frac{\mu_0}{4\pi} \int_{\Omega_{s}} \frac{\Js(\rv')}{\left| \rv - \rv'  \right|} \mathrm{d} \vec{s} 
		\label{eq:BiotSavart}
\end{equation}	
where $\mu_0$ is the permeability in $\Vo$, $\Omega_{s}$ is the domain where the source current $\Js$ is present, in our case the source coils, and $\rv$ the evaluation point. For source coil windings, it is justified to model the source regions as one‑dimensional lines by assuming that the current density is constant along the small cross-section~\cite{Schops_2013aa}. Note that the integration of these line currents poses no difficulty thanks to Biot-Savart's law \eqref{eq:BiotSavart}. Their paths may be chosen arbitrarily, reflecting the flexibility of practical winding layouts. Recall that
\begin{equation}
    \gamma_\mathrm{D}(\vec{u}) = \vec{u} \times \vec{n},
    \qquad
    \gamma_\mathrm{N}(\vec{u}) = \mu _{0}^{-1} (\CURL \vec{u}) \times \vec{n} ,
\end{equation}
denote the tangential Dirichlet and Neumann traces commonly used in curl‑conforming formulations, where \(\vec{n}\) is the outward unit normal~\cite{Monk_2003aa}.
The image field $\Am$ is computed as: \newline
Find $\Am \in \{\vec{v} \in  \HSPACE{}{}{\CURL}{\Vo}: \gamma_{\mathrm{D}}(\vec{v})  = -\gamma_{\mathrm{D}}(\As)  \text{ on } \Gio, \gamma_{\mathrm{D}}(\vec{v})  = 0 \text{ on } \partial V\}$ such that
\begin{align}
	\left(\mu _{0}^{-1} \CURL \Am, \, \CURL \Am'\right)_{\Vo} & = 0
    \label{eq:image_problem}
\end{align}
for all $\Am' \in \HSPACE{}{0}{\CURL}{\Vo}$, where $\Am'$ is a corresponding test function and $\gamma_{\mathrm{D}}(\As) \in \HSPACE{-1/2}{}{\CURL}{\Gio}$, the corresponding trace space. The surface restriction of the source field follows from an exterior representation formula together with the associated jump and continuity conditions that ensure equivalence between the exterior and interior formulations outside the coil region. In this framework, we define the virtual surface current density $ \vec{K}_{\mathrm{g}} = \gamma_\mathrm{N} (\Am + \As)$ on $\Gio$, which represents the Neumann‑type tangential trace of the superposition of image and source field on the interface $\Gio$. To compute $\gamma_\mathrm{N} (\As)$, it is either possible to compute the magnetic field $\vec{B}$ via Biot-Savart law~\cite[Sec.~5.2]{Jackson_1998aa}
\begin{equation}
    \vec{B}_\mathrm{s}(\rv) = \frac{\mu_0}{4\pi} \int_{\Omega_{s}} \frac{\Js(\rv) \times (\rv-\rv')}{\left| \rv - \rv'  \right|^3} \mathrm{d} \vec{s},
\end{equation} and take its tangential components, or to compute the $L^2$-projection of \eqref{eq:BiotSavart} onto the discrete curl-conforming space and then take its Neumann trace, \cite{DAngelo_2024ab}. The reaction field $\Ag$ must satisfy the corresponding jump condition across $\Gio$, in particular $\gamma_N^+ (\Ag) - \gamma_N^- (\Ag) = \Kg$. The weak formulation to compute the reaction field $\Ag(t)$ therefore reads:\newline Find $\Ag \in \HSPACE{}{0}{\CURL}{V}$ subject to
\begin{align}
    \left( \mu^{-1} \CURL \Ag , \CURL \Ag'\right)_V  &- \mathrm{j} \omega \left(\sigma \Ag, \Ag'\right)_V
     \label{eq:red_problem}     \\
    &\hspace{3em} = \langle \vec{K}_{\mathrm{g}}, \gamma_\mathrm{D}(\Ag') \rangle_{\Gio},
    \nonumber   
\end{align}
for all $\Ag' \in \HSPACE{}{0}{\CURL}{V}$, where $\gamma_\mathrm{D} (\Ag') \in  \HSPACE{-1/2}{}{\CURL}{\Gio}$ denotes the tangential trace on $\Gio$ and $ \vec{K}_{\mathrm{g}} \in  \HSPACE{-1/2}{}{\DIV}{\Gio}$. The duality term on the boundary is well‑defined because $\HSPACE{-1/2}{}{\DIV}{\Gio}$ is the dual of $ \HSPACE{1/2}{}{\CURL}{\Gio}$, so the pairing $\langle \vec{K}_{\mathrm{g}}, \gamma_\mathrm{D}(\Ag') \rangle_{\Gio}$ fits naturally within the trace framework of $\HSpace{}{}{\CURL}$-conforming spaces.
The location of the integration surface $\Gamma$ is not unique and may be chosen freely, as long as it separates the linear source region from the conducting domain. While placing the interface close to the conductors is admissible, care must be taken to avoid proximity to the source coils, as this may deteriorate the accuracy of the Biot–Savart evaluation due to (near-) singular integrals. The placement can also be guided by the region of interest, since the field inside $\Gamma$ is represented solely by $\Ag$, whereas outside additional contributions must be evaluated.There are several advantages to the proposed surface-reduced formulation. First, the source field to be reused in parametric studies as long as geometry changes occur only inside $\Vi$. Second, in the presence of nonlinear materials, the decomposition into linear and nonlinear field components implies that the linear source field remains unchanged throughout the Newton iteration. In both cases, the source field needs to be computed only once on the interface, whereas the reaction field alone is updated for the parametric instances or during the nonlinear solve, respectively. Moreover, because the source field is evaluated exclusively outside the coil region, all singularities inherent to the Biot-Savart kernel are avoided, and only smooth integrals are encountered on the integration surface -- significantly simplifying numerical quadrature. Finally, the surface‑reduced formulation offers a substantial runtime advantage. In the original RMVP, the Biot-Savart kernel must be evaluated at every quadrature point in the 3D finite-element mesh, leading to  $\mathcal{O}(nm^3)$ work for a tensor-product discretization with $m^3$ quadrature points in the volume and $n$ quadrature points along the coil. Implementation of the FMM would allow to reduce the cost to at most $\mathcal{O}(n+m^3)$~\cite{Carrier_1988aa}. In contrast, the interface RMVP requires Biot-Savart evaluations only on the separating interface, which contains merely $\mathcal{O}(m^2)$ quadrature points, yielding a total cost of $\mathcal{O}(nm^2)$ without any multi-pole machinery. For $n \approx m$, this yields the same runtime complexity as the original approach with FMM.  Additional evaluations are required only at points outside the interface, where the solution is sought. The price for this reduction is a single auxiliary magnetostatic solve inside the air region, the image problem, which is cheap compared to the large number of Biot-Savart evaluations otherwise required, which become increasingly costly for coils with many turns.
\subsection{IGA Discretization}
\label{sec:discretization}
Our geometry representation is inspired by CAD. It is given by a map from the reference domain $\hat{\Omega} = \left[0,1\right]^3$ to a physical domain $\Omega \subset \mathbb{R}^3$. The standard tools for this representation are B-splines and NURBS. Given a knot vector $\Xi 
\subset \left[0,1\right]$, the basis of univariate B-splines $\hat B_i^p$ of degree $p$ can be defined using the Cox-de-Boor recursion formula. From those one derives the NURBS basis functions $\hat{N}_i^p$ \cite{Cohen_2001aa}.
Curves are then described as linear combination of these functions, surfaces and volumes are created using tensor products~\cite{Cohen_2001aa}. This yields a higher-order, smooth representation of the domain, in particular it is able to exactly represent many curved geometries, include all conic section~\cite{Hughes_2005aa}. Note that, in most practical applications, the geometry cannot be parametrized using a single map from the reference to the physical domain. This is also the case for our setup~(see. \autoref{fig:problem_sketch}), which features at least two regions with different materials. In such cases, a multi-patch parametrization is used, where the physical domain is decomposed into a collection of subdomains, each with a corresponding map, to be appropriately combined~\cite{Buffa_2015aa}. The surface $\Gio$ will be placed exactly at a patch interface and is therefore mesh-aligned, the coil is not resolved in this mesh. It is modeled as a closed one-dimensional spline curve, again allowing for smooth representation of curved paths. For the numerical evaluation of \eqref{eq:BiotSavart}, we employ the trapezoidal rule, which exhibits exponential convergence for smooth periodic integrands on closed curves~\cite{Trefethen_2014aa}. Since physical coils form current loops, the assumption of closed curves is naturally satisfied. The trapezoidal rule requires equidistant quadrature points in the physical domain. Due to the non-linear mapping from the reference to the physical space, these points must be determined by solving a nonlinear system once at a preprocessing stage.
Finally, IGA is used to discretize the subproblems \eqref{eq:image_problem} and \eqref{eq:red_problem}. Following the standard Ritz-Galerkin approach the B-splines are used to span the ansatz and test function spaces. The magnetic vector potential is approximated by
\begin{equation}
\vec{A}(\textbf{x})
\approx 
\sum\nolimits_{i=1}^N \vec{B}_i^p(\textbf{x})\; u_i
\label{eq:AnsatzTestFunctions}
\end{equation}
where  $u_i$ are the unknown coefficients and $\vec{B}_i^p(\textbf{x})$ are the spline basis functions of order $p$. We employ the spline-based $\HSPACE{}{}{\CURL}{V}$-conforming space $S_p^1(V)$ of degree $p$, which is the analog to classical Nédélec edge elements, see~\cite{Buffa_2019ac}. In particular, we denote the volumetric space and the trace space of the Dirichlet and Neumann trace as
\begin{align*}
    S_p^1(V) \subset \HSPACE{}{}{\CURL}{V} \\
    \gamma_\mathrm{D} :  S_p^1(V) \rightarrow S_p^1(\Gamma) \\
    \gamma_\mathrm{N} : S_p^1(V) \rightarrow S_p^{1\ast}(\Gamma),
\end{align*} respectively. Tree-cotree gauging is applied to ensure a uniquely solvable system~\cite{Albanese_1988aa, Kapidani_2022aa}. 
\subsection{Error estimates}
\label{sec:error_estimates}
For the standard convergence theory of isogeometric analysis to apply, the underlying exact solution has to be sufficiently regular, as detailed in \cite{Buffa_2010aa}. In the variational formulation~\eqref{eq:red_problem}, the right-hand side consists of a surface distribution, represented by the Maxwell-Neumann trace. This is a considerably milder singularity structure compared to, for instance, Dirac delta sources. As a consequence, the continuous solution $\Ag$ remains piecewise smooth provided that the interface $\Gio$ is Lipschitz-continuous. More precisely, although $\Ag \in \HSPACE{}{}{\CURL}{V}$ globally, it enjoys higher Sobolev regularity within each of the subdomains seperated by $\Gio$. In particular, it has local regularity $\HSpace{k+1}{}{V_\mathrm{i}}$ for $i \in \{\mathrm{int},\mathrm{ext}\}$, as established in the classical interface-regularity results \cite{Bramble_1996aa}. Such piecewise smoothness is sufficient to achieve optimal convergence of curl-conforming edge elements on shape regular meshes, because the interpolation operators underlying the error estimates act locally. Under these assumptions, we therefore expect the discrete solution $\Agh$ with mesh size $h$ to satisfy the optimal convergence rates \cite[Th.~5.4]{Buffa_2015aa}
\begin{align}
\|\Ag - \Agh\|_{\HSPACE{}{}{\CURL}{V}} &= \mathcal{O}(h^{p}),
\label{eq:rateHcurl}
\end{align}
measured in the \HSpace{}{}{\CURL}-norm and, equivalently, in the $\HSpace{}{}{\CURL}$-seminorm
\begin{align}
    |\Ag - \Agh|_{\HSPACE{}{}{\CURL}{V}} &:= \|\CURL(\Ag - \Agh)\|_{L^{2}(V)} \nonumber \\ &= \mathcal{O}(h^{p}). \label{eq:rateHcurlsemi}
\end{align}
The $\HSpace{}{}{\CURL}$-seminorm is gauge-invariant and therefore provides a robust and physically meaningful measure of the discretization error. For this reason, all numerical results reported in this work focus on $|\Ag - \Agh|_{\HSPACE{}{}{\CURL}{V}}$, which directly reflects the accuracy of the magnetic field $\vec{B}=\CURL \Ag$.
The convergence rates  stated in (\ref{eq:rateHcurl}--\ref {eq:rateHcurlsemi}) are attainable, as long as the computational mesh is aligned with the interface $\Gio$, a requirement well documented in the theory of elliptic interface problems\,\cite{Babuska_1970aa,Barrett_1987aa}. In essence, misalignment would deteriorate the local approximation properties, whereas alignment ensures that the piecewise smoothness of the solution is properly captured by the discrete space. Furthermore, achieving this convergence rates in practice requires that the auxiliary fields $\As$ and $\Am$ entering the formulation are computed with sufficiently high accuracy. In particular, the kernel evaluation~\eqref{eq:BiotSavart} must rely on a fast converging quadrature rule to not deteriorate the approximation error. Likewise, the discretization used for solving the image problem~\eqref{eq:image_problem} must employ spline spaces of degree and regularity at least matching those used for $\Ag$, so as not to create a bottleneck in the approximation pipeline. A rigorous proof is beyond the scope of this work. Nevertheless, the numerical experiments presented in the next section will confirm the theoretical expectations and illustrate the practical effectiveness of the proposed isogeometric discretization. 
\section{Numerical Studies}

\label{sec:results}
First, we verify the implementation of the proposed method against an analytical reference solution and demonstrate the expected convergence rates. We then examine in detail the components required to actually attain these rates, in particular the quadrature strategy for evaluating the kernel and the choice of discrete spaces for computing $\As$ and $\Am$. Finally, we apply the method to a more complex configuration, compare the results with a commercial simulation tool, and confirm the effectiveness of the proposed approach.  All computations rely on \textit{GeoPDEs} \cite{Vazquez_2016aa} and we provide our implementations in \cite{Backmeyer_2026aa}.
\subsection{Verification}
\begin{figure}[t]
	\centering
	\begin{tikzpicture}[scale=0.75]

\begin{axis}[
	tudalineplot,
	height = 7cm,
    width = 9cm,
	axis line style = {line width=2pt},
	at={(0.758in,0.481in)},
	scale only axis,
	xmode=log,
	xmin=10,
	xmax=87,
	ymode=log,
	ymin=0.0013,
	ymax=1,
	ylabel style={font=\color{white!15!black}},
	ylabel={$|\Ag - \Agh|_{\HSPACE{}{}{\CURL}{\Vi}}$},
	xlabel={$\sqrt[3]{\text{\#dofs}}$},
	axis background/.style={fill=white},
	axis x line*=bottom,
	axis y line*=left,
	grid=major,
	legend style={at={(0.0,0.0)},anchor=south west},
	legend columns=2,
	]
\addplot [color=TUDa-1b, mark=*, line width=1.5pt, mark size=3pt]
table[row sep=crcr]{
12.0690463987073	0.273253890900853\\
17.8631598770806	0.192056916605968\\
23.6595451252288	0.149651342528766\\
29.4568938560715	0.122911062600901\\
35.2547400697818	0.104817586677744\\
41.0528764533186	0.0915864441057624\\
46.8511968036355	0.0814630364545838\\
52.6496410989515	0.073436198342954\\
58.4481728625957	0.0669251534266559\\
64.2467686505402	0.0615234284684133\\
70.0454127088693	0.0569488053441396\\
75.8440940599367	0.0530337953462229\\
81.6428048199676	0.0496432698179854\\
};
\addlegendentry{$p=1$}
\addplot [color=TUDa-1b, mark=none, line width= 1.5pt, dotted, mark size=3pt]
table[row sep=crcr]{
6.28276130478279	0.529533596869315\\
12.0690463987073	0.27565833141128\\
17.8631598770806	0.186245502749019\\
23.6595451252288	0.140616954991471\\
29.4568938560715	0.112942430666608\\
35.2547400697818	0.0943683937367793\\
41.0528764533186	0.0810401969220366\\
46.8511968036355	0.0710106340706076\\
52.6496410989515	0.0631900450326002\\
58.4481728625957	0.0569210811741572\\
64.2467686505402	0.0517836657293893\\
70.0454127088693	0.047496803335584\\
75.8440940599367	0.0438654219979717\\
81.6428048199676	0.0407498639877608\\
};
\addlegendentry{$\mathcal O(h^1)$}
\addplot [color=TUDa-9c, mark=diamond*, mark size=3pt, line width=1.5pt]
table[row sep=crcr]{
11.6789892476198	0.254789969289569\\
17.0677801465357	0.103375852784154\\
22.8700931077745	0.0543843268409115\\
28.6708725275452	0.0356898642617871\\
34.4709550405101	0.0263825624364454\\
40.2706643887459	0.0201434267415163\\
46.0701513119509	0.0160785217178729\\
51.86949523046	0.0132352885556199\\
57.6687418738576	0.0111946650579824\\
63.4679193979439	0.00970715755222394\\
69.2670460715456	0.00861987321909473\\
75.0661342589156	0.00783105052767735\\
80.8651926232716	0.00726508317312339\\
86.6642274132671	0.00686490693616334\\
};
\addlegendentry{$p=2$}
\addplot [line width = 1.5pt,TUDa-9c,mark = none,dotted,mark size=3pt] 
table[row sep=crcr]{
11.6789892476198	0.254789969289569\\
17.0677801465357	0.119299538509934\\
22.8700931077745	0.0664441849251536\\
28.6708725275452	0.0422776680008771\\
34.4709550405101	0.0292473290815952\\
40.2706643887459	0.0214296586989274\\
46.0701513119509	0.0163739540149536\\
51.86949523046	0.0129172123855081\\
57.6687418738576	0.0104498941203228\\
63.4679193979439	0.00862748713982672\\
69.2670460715456	0.00724335063573416\\
75.0661342589156	0.00616743727416938\\
80.8651926232716	0.00531458788130643\\
86.6642274132671	0.00462714494691622\\
};
\addlegendentry{$\mathcal O(h^2)$}
\addplot [color=TUDa-3c, line width=1.5pt,mark=square*, mark size=3pt]
table[row sep=crcr]{
16.1904328004549	0.105095429973991\\
22.0220165270767	0.0610866972385944\\
27.8392145372266	0.0249808197247817\\
33.6498080262115	0.0146794386118852\\
39.4568240511079	0.00925011511013523\\
45.2616851683682	0.00728621557992437\\
51.0651480648635	0.00638037240682757\\
56.8676523132809	0.0060422777804805\\
62.6694707260955	0.00590512300017262\\
68.4707815785799	0.00584215629283683\\
};
\addlegendentry{$p=3$}
\addplot [color=TUDa-3c, line width= 1.5pt, dotted, mark=none, mark size = 3pt]
table[row sep=crcr]{
16.1904328004549	0.105095429973991\\
22.0220165270767	0.0417626408997768\\
27.8392145372266	0.020672274972637\\
33.6498080262115	0.0117060785473103\\
39.4568240511079	0.00726093971494458\\
45.2616851683682	0.00481024335457504\\
51.0651480648635	0.00334954194059489\\
56.8676523132809	0.00242528794498067\\
62.6694707260955	0.00181213740923413\\
68.4707815785799	0.00138945136822016\\
};
\addlegendentry{$\mathcal O(h^3)$}
\end{axis}
\end{tikzpicture}%
	\caption{Convergence results for degree $p$.}
	\label{fig:cvg_study}
\end{figure}
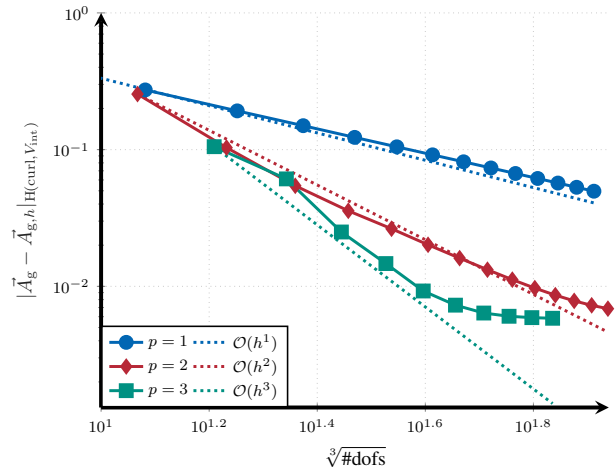

To verify the correctness of the proposed implementation, we first compare numerical results against an existing analytical solution for a canonical eddy current configuration described in \cite{Bowler_2005aa}. The setup consists of a conducting cylinder of radius $r=\SI{12}{\milli \meter}$ and height $h=\SI{60}{\milli\meter}$, centered at the origin and surrounded by air. The cylinder has an electrical conductivity of $\sigma=\SI{35e6}{\siemens \per \meter}$ and is excited by a circular coil of radius $r_\mathrm{coil} = \SI{25}{\milli\meter}$ carrying a sinusoidal current with amplitude $I = \SI{320}{\ampere} $ at frequency $f=\SI{200}{\hertz}$. In the convergence study shown in \autoref{fig:cvg_study}, we can see that the solution converges with the optimal convergence rates up to a threshold, which is the best achieved accuracy with the superposition of Bessel functions in the reference solution. The same rates are also observe if the error is computed in $V$, where in $\Vo$ the three fields $\As$, $\Am$ and $\Ag$ act in superposition, if we exclude the singularity at the source coils from the computation domain. 
\subsection{Quadrature of the source field}
\label{sec:quad}
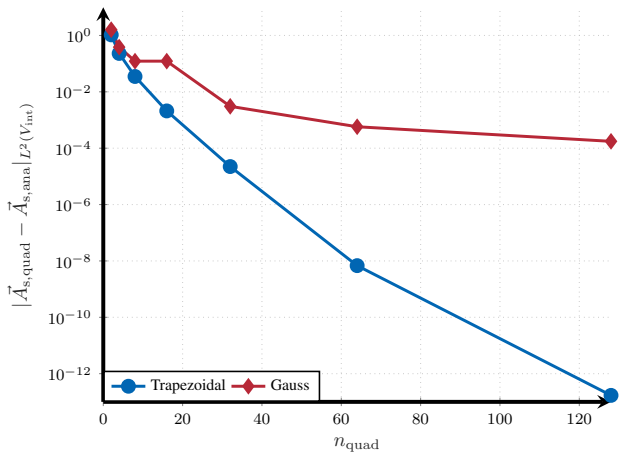
\begin{figure}[t]
	\centering
	\begin{tikzpicture}[scale=0.75]

\begin{axis}[
    tudalineplot,
    height=7cm,
    width=9cm,
    axis line style = {line width=2pt},
	at={(0.758in,0.481in)},
	scale only axis,
    xmin=0,
    xmax=128,
    xlabel style={font=\color{white!15!black}},
    xlabel={$n_\mathrm {quad}$},
    ymode=log,
    ymin=1e-13,
    ymax=10,
    ylabel style={font=\color{white!15!black}},
    ylabel={$|\vec{A}_\mathrm{s,quad} - \vec{A}_\mathrm{s,ana}|_{L^2(\Vi)}$},
    axis background/.style={fill=white},
    legend style={legend cell align=left, align=left},
    axis x line*=bottom,
	axis y line*=left,
	grid=major,
	legend style={at={(0.0,0.0)},anchor=south west},
	legend columns=2,
]
\addplot [color=TUDa-1b, mark=*, line width=1.5pt, mark size=3pt]
  table[row sep=crcr]{
2	1.05\\
4	0.23\\
8	0.035\\
16	0.0021\\
32	2.23e-05\\
64  6.81e-09\\
128 1.7e-13\\
};
\addlegendentry{Trapezoidal}

\addplot [color=TUDa-9c, mark=diamond*, mark size=3pt, line width=1.5pt]
  table[row sep=crcr]{
2	1.60146126011282\\
4	0.38507731730916\\
8	0.122501288458325\\
16	0.01183489530508\\
32	0.003057124213876\\
64  5.766470505136153e-04\\
128 1.751072459518974e-04\\
};
\addlegendentry{Gauss}

\end{axis}
\end{tikzpicture}%
	\caption{Convergence of the error in the source field for the different quadrature rules.}
	\label{fig:quad}
\end{figure}
In the previous test case, the coil source was simply a circle. In that case, numerical evaluation of the integral appearing in~\eqref{eq:BiotSavart} would not have been necessary, since for such a simple setup there exist closed-form expression for the magnetic vector potential $\vec{A}_\mathrm{s,ana}$ ~\cite[Sec.~5.5]{Jackson_1998aa}. We use this analytical formula to assess the quality of our quadrature rule; as stated in \autoref{sec:error_estimates} it is essential to evaluate the integral with high accuracy. \autoref{fig:quad} shows the expected spectral (exponential) convergence of the quadrature error applying trapezoidal rule to the Biot–Savart integral over a closed curve, and contrasts this with the merely algebraic convergence of Gaussian quadrature when refining the number of quadrature points. When the kernel is not approximated accurately enough, the resulting quadrature error ultimately limits the achievable solution accuracy. This underscores the need for a rapidly converging quadrature rule. In all computation, we set $n_\mathrm{quad}=64$, which proved sufficient to ensure that quadrature error does not limit the overall convergence. For practical simulations, one should implement an adaptive procedure, e.g. by measuring the change in $\Kg$ to determine the number of required quadrature points.
\subsection{Trace‑space requirements for high‑order convergence}
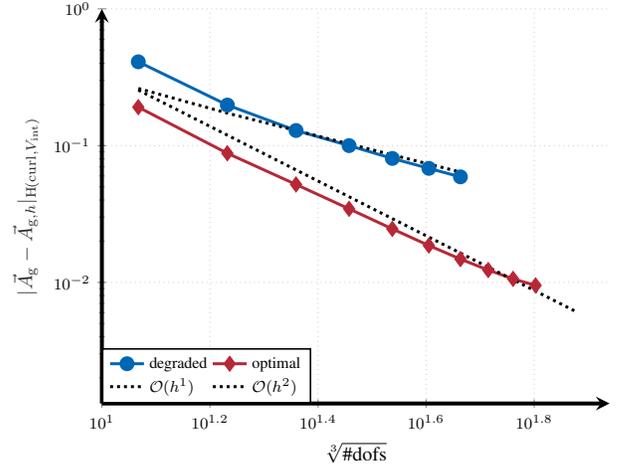
\begin{figure}[t]
	\centering
	\begin{tikzpicture}[scale=0.75]

\begin{axis}[
	tudalineplot,
	height = 7cm,
    width = 9cm,
	axis line style = {line width=2pt},
	at={(0.758in,0.481in)},
	scale only axis,
	xmode=log,
	xmin=10,
	xmax=87,
	ymode=log,
	ymin=0.0013,
	ymax=1,
	ylabel style={font=\color{white!15!black}},
	ylabel={$|\Ag - \Agh|_{\HSPACE{}{}{\CURL}{\Vi}}$},
	xlabel={$\sqrt[3]{\text{\#dofs}}$},
	axis background/.style={fill=white},
	axis x line*=bottom,
	axis y line*=left,
	grid=major,
	legend style={at={(0.0,0.0)},anchor=south west},
	legend columns=2,
	]
\addplot [color=TUDa-1b, mark=*, mark size=3pt, line width=1.5pt]
table[row sep=crcr]{
11.6789892476198	0.4103\\
17.0677801465357	0.1981\\
22.8700931077745	0.1290\\
28.6708725275452	0.1\\
34.4709550405101	0.0806\\
40.2706643887459	0.0682\\
46.0701513119509	0.0591\\
};
\addlegendentry{degraded}

\addplot [color=TUDa-9c, mark=diamond*, mark size=3pt, line width=1.5pt]
table[row sep=crcr]{
11.6789892476198	0.1911\\
17.0677801465357	0.0879\\
22.8700931077745	0.0520\\
28.6708725275452	0.0345\\
34.4709550405101	0.0246\\
40.2706643887459	0.0186\\
46.0701513119509	0.0148\\
51.86949523046	    0.0123\\
57.6687418738576	0.0106\\
63.4679193979439	0.0095\\
};
\addlegendentry{optimal}

\addplot [line width = 1.5pt,black,mark = none,dotted,mark size=3pt] 
table[row sep=crcr]{
11.6789892476198	0.2620\\
17.0677801465357	0.1729\\
22.8700931077745	0.1290\\
28.6708725275452	0.1029\\
34.4709550405101	0.0856\\
40.2706643887459	0.0733\\
46.0701513119509	0.0641\\
};
\addlegendentry{$\mathcal O(h^1)$}

\addplot [line width = 1.5pt,black,mark = none,dotted,mark size=3pt] 
table[row sep=crcr]{
11.6789892476198	0.254789969289569\\
17.0677801465357	0.119299538509934\\
22.8700931077745	0.0664441849251536\\
28.6708725275452	0.0422776680008771\\
34.4709550405101	0.0292473290815952\\
40.2706643887459	0.0214296586989274\\
46.0701513119509	0.0163739540149536\\
51.86949523046	0.0129172123855081\\
57.6687418738576	0.0104498941203228\\
63.4679193979439	0.00862748713982672\\
69.2670460715456	0.00724335063573416\\
75.0661342589156	0.00616743727416938\\
};
\addlegendentry{$\mathcal O(h^2)$}

\end{axis}
\end{tikzpicture}%
	\caption{Convergence of the reduced field in the $\HSpace{}{}{\CURL}$–seminorm under trace space mismatch. The bulk for $\Ag$ is $S_p^{1}(V)$ (here: $p=2$), while the interface current density $\Kg$ is approximated either in the matching trace space $S_{p}^{1\ast}(\Gio)$ (optimal) or in the degraded space $S_{p-1}^{1\ast}(\Gio)$.}
	\label{fig:cvg_red}
\end{figure}

A central feature of the proposed method is the interface term arising from the exterior‑representation principle.
As discussed in \autoref{sec:error_estimates}, the discrete trace space used for $\Kgh = \gamma_\mathrm{N}(\As+\Am)$ must match the trace space of $\Agh$ in both degree and regularity. To demonstrate this numerically, we perform an experiment in which we deliberately reduce the polynomial degree of the trace space used for $\Kgh$ by one order, and then study the resulting convergence rates of  $\|\Ag - \Agh\|_{\HSPACE{}{}{\CURL}{\Vi}}$. The results in ~\autoref{fig:cvg_red} compare two cases, both with $\Agh$ is represented in $S_p^1(V)$. When the surface current density is approximated in the matching trace space $S_p^{1\ast}(\Gio)$ (the Neumann trace of $S_p^1(V)$, we recover the optimal rate. In contrast, when $\Kgh \in  S_{p-1}^{1\ast}(\Gio)$, the convergence rate degrades accordingly (approximately $\mathcal{O}(h^{p-1})$).
This confirms that the attainable convergence order of the reduced field $\Agh$ is directly limited by the approximation order of the trace space used for $\Kg$. 
Note that the error is evaluated only inside $\Vi$, where $\Agh$ represents the full field according to~\eqref{eq:field_superpos}. As a result, any inaccuracy in approximating the source or image field propagates into $\Vi$ via the interface term. Consequently, the overall convergence is restricted by the polynomial degree of the employed trace space.
\subsection{Flexible coil path}
\begin{figure}[t]
    \centering
	\includegraphics[width=0.75\linewidth]{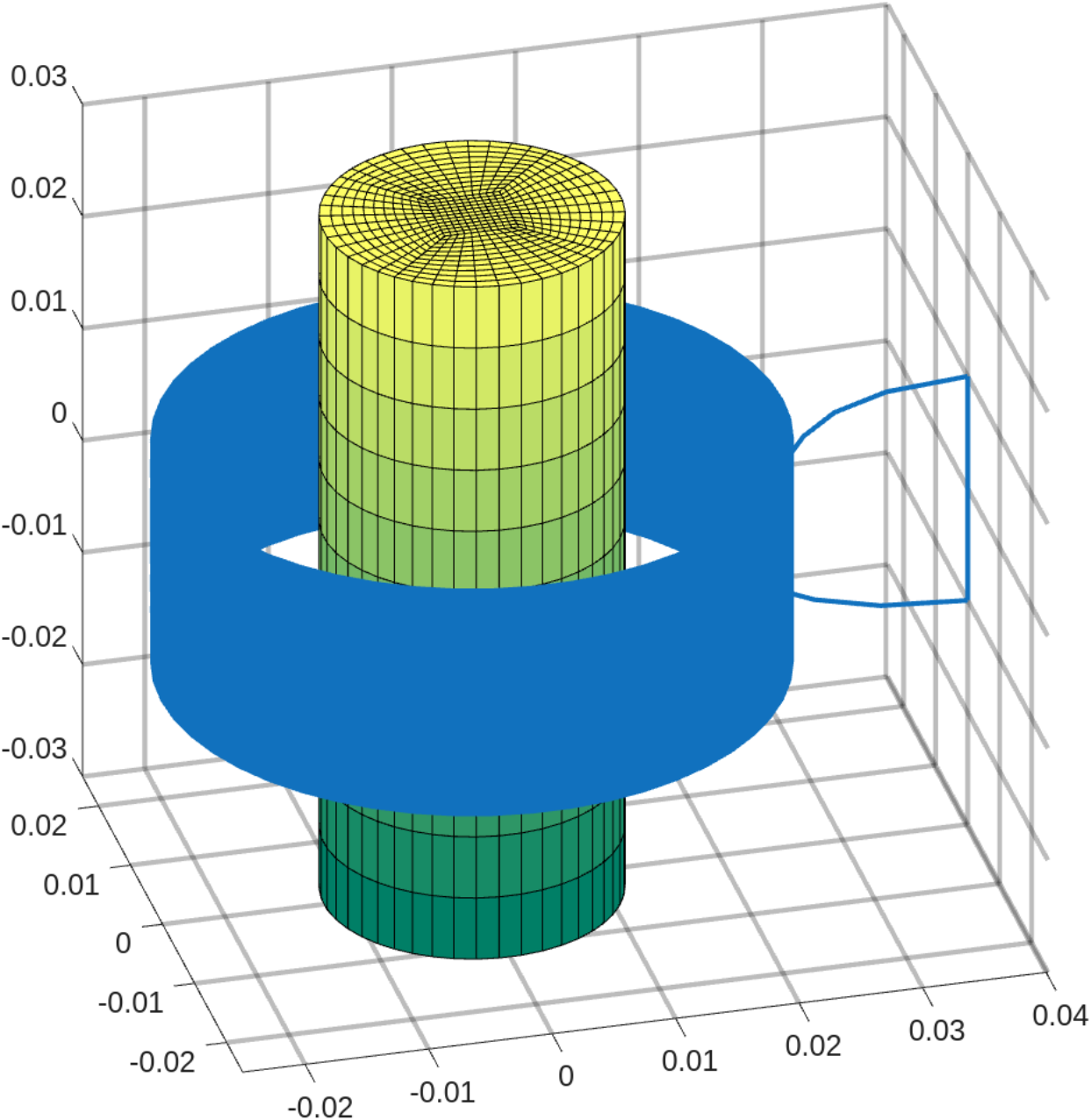}
	\caption{Problem setup with helicoidal coil of 200 turns wound around cylindrical conducting specimen.}
	\label{fig:complex_coil}
\end{figure}
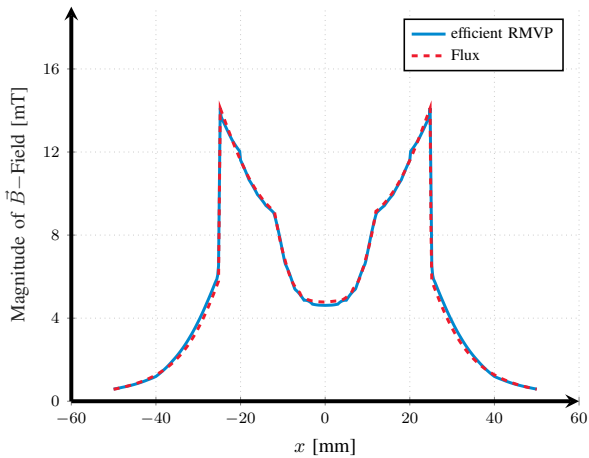
\begin{figure}[t]
	\centering
	\begin{tikzpicture}[scale=0.75]
  \begin{axis}[
	tudalineplot,
    height = 7cm,
    width = 9cm,
	axis line style = {line width=2pt},
	at={(0.758in,0.481in)},
	scale only axis,
    ylabel style={font=\color{white!15!black}},   
    xmin=-60, xmax=60,
    ymin=0,   ymax=19,
    xlabel={$x$ [mm]},
    ylabel={Magnitude of $\vec{B}-$Field [mT]},
    axis x line*=bottom,
	axis y line*=left,
	grid=major,
    xtick distance=20,
    ytick distance=4,
    legend pos=north east,
  ]
    \addplot[TUDa-2b,line width=1.5pt,]
    table[col sep=comma,
      x expr=\thisrowno{0} * 1000,
      y expr=\thisrowno{1} * 1000] {data/Ypath_values.csv};
    \addlegendentry{efficient RMVP}

    \addplot[TUDa-9b,dashed,line width=1.5pt,]
    table[col sep=comma,
      x expr=\thisrowno{0} * 1000,
      y expr=\thisrowno{1} * 1000]{data/FLUX_PathY_cleaned.csv};
    \addlegendentry{Flux}
  \end{axis}
\end{tikzpicture}
	\caption{Magnitude of the magnetic flux density $B$ along a horizontal axis compared against \textit{Flux}.}
	\label{fig:flux_comparison}
\end{figure}

The methodology is subsequently extended to a more complex configuration shown in~\autoref{fig:complex_coil}: a helicoidal coil comprising 200 turns wound around the same conducting cylindrical specimen. The computed results are benchmarked against those obtained using the commercial electromagnetic simulation software \textit{Flux}\footnote{\url{https://altair.com/flux}}. In \autoref{fig:flux_comparison}, the magnetic flux density is compared along the horizontal axis. The relative difference of the magnetic energy inside the conducting cylinder computed by the two implementations is below $2\,\%$. This error is to be expected due to different meshes, basis functions and implementations. By applying the interface RMVP, the number of kernel evaluations (KE) reduces by over $85\%$ for a discretization with approx. $18000$ degrees of freedom ($9792$ (interface) vs. $79488$ (original)) to obtain the field solution in $\Vi$. To obtain the field solution in the full computational domain $V$, it could still be reduced by over $50\%$. A naive implementation leads to runtimes of $\SI{34}{\second}$ -- including $\SI{10}{\second}$ dedicated to the KE -- and $\SI{100}{\second}$ ($\SI{85}{\second}$ for KE), for the interface and original RMVP approach respectively. When optimizing the original version by vectorization we can reduce the runtime to $\SI{72}{\second}$ ($\SI{43}{\second}$ for KE). A similar optimization has not yet been implemented in the interface method.
\section{Conclusion}
\label{sec:conclusion}
In this work, we revisited the interface RMVP formulation introduced in \cite{DAngelo_2024ab}, in which the evaluation of the Biot-Savart kernel is restricted to a separating surface in order to reduce computational cost. While \cite{DAngelo_2024ab} demonstrated the efficiency of this approach for magnetostatic problems with low-order discretizations, the reported convergence behavior remained suboptimal.
The key novelty of the present contribution is the detailed analysis of the convergence behavior. We show that, contrary to earlier observations, the interface formulation does not limit the convergence order. Instead, optimal convergence rates are achieved, even for higher-order discretizations. This enables the method to exploit the smoothness and geometric exactness inherent to spline spaces, yielding in higher-order accuracy while simplifying the representation of complex coil geometries in multipatch geometries. Numerical experiments with spline‑based discretizations support this conclusion through validation against an analytical reference solution. The same convergence behavior can be expected for conventional high‑order finite element methods. In addition, the results highlight that suitable trace spaces are essential for the accurate approximation of the interface term, and, consequently, for achieving the expected convergence rates.
Furthermore, we present an effcient quadrature strategy for the Biot-Savart kernel based on the trapezoidal rule, which achieves exponential convergence. Finally, this work operates in the MQS regime and utilizes arbitrary winding paths. Overall, we have demonstrated that, even in the higher-order setting, the interface RMVP formulation significantly reduces the number of kernel evaluations, while still attaining optimal convergence rates. 
\section*{Acknowledgment}
The authors thank Herbert Egger and Herbert De Gersem for the fruitful discussions. 
This work was supported by the Graduate School CE at TU Darmstadt, 
the IRGA program of Université Grenoble Alpes, the German Research Foundation (DFG) and the Austrian Science Fund (FWF) through the Collaborative Research Centre/Transregio TRR 361/F90 CREATOR (DFG Project-ID 492661287).

 \end{document}